\documentclass[12pt]{article}
\usepackage{amsfonts, amsmath, amsthm, amssymb}
\usepackage[all]{xy}

\newtheorem{theorem}{Theorem}[section]
\newtheorem{lemma}[theorem]{Lemma}

\newtheorem{prop}[theorem]{Proposition}

\def\NN{\mathbb{N}}
\def\PP{\mathbb{P}}

\def\ZZ{\mathbb{Z}}
\def\calE{\mathcal{E}}
\def\calF{\mathcal{F}}

\def\calO{\mathcal{O}}
\def\calR{\mathcal{R}}

\def\alg{\mathrm{alg}}
\def\an{\mathrm{an}}
\def\con{\mathrm{con}}

\def\Gancon{\Gamma_{\an,\con}}
\def\Gcon{\Gamma_{\con}}

\def\be{\mathbf{e}}

\def\bv{\mathbf{v}}
\def\bw{\mathbf{w}}

\def\gotha{\mathfrak{a}}

\DeclareMathOperator{\divi}{div}
\DeclareMathOperator{\Frac}{Frac}

\DeclareMathOperator{\inte}{int}

\DeclareMathOperator{\Span}{SatSpan}
\DeclareMathOperator{\Spec}{Spec}

\newcounter{fixmectr}

\begin{document}

\title{Semistable reduction for overconvergent $F$-isocrystals
on a curve}
\author{Kiran S. Kedlaya \\
Department of Mathematics \\ University of California, Berkeley \\
Berkeley, CA 94720}
\date{March 5, 2003}

\maketitle

\begin{abstract}

Let $X$ be a smooth affine curve over a field $k$ of characteristic $p>0$
and $\calE$ an overconvergent $F^a$-isocrystal on $X$ for some positive
integer $a$. We prove that after replacing $k$ by some finite purely
inseparable extension, there exists a finite separable
morphism $X' \to X$,
the pullback of $\calE$ along which extends to
a log-$F^a$-isocrystal on a smooth compactification of $X'$.
This resolves a weak form of
the global version of a conjecture of Crew; the proof uses the local version
of the conjecture, established (separately)
by Andr\'e, Mebkhout and the author.
\end{abstract}

\section{Introduction}

The purpose of this paper is to prove the following
semistable reduction theorem
for overconvergent $F$-isocrystals on a curve,
answering a conjecture of de Jong \cite[Section~5]{bib:dej2}, a
reformulation of a conjecture of Crew (of which more about below).
\begin{theorem}[Semistable reduction] \label{thm:main}
Let $X$ be a smooth, geometrically connected curve over a field $k$ 
of characteristic $p>0$, let $K$ be a finite, totally ramified extension
of $\Frac(W(k))$ admitting a lift of the $p$-power Frobenius on $k$,
and let $\calE$ be an
overconvergent $F^a$-isocrystal on $X$ with respect to $K$ 
for some $a \in \NN$.
Then after replacing $k$ by a suitable finite purely inseparable
extension (depending on $\calE$), there exist
a finite generically \'etale 
morphism $f: X_1 \to X$, a smooth compactification
$j: X_1 \hookrightarrow \overline{X_1}$ of $X_1$, and
a log-$F^a$-isocrystal $\calF$ on $(\overline{X_1},\overline{X_1} \setminus
X_1)/K$
such that $j^{*}\calF \cong f^* \calE$.
\end{theorem}

Our basic approach is to use 
the quasi-unipotence theorem
($p$-adic local monodromy theorem) for $F^a$-isocrystals,
plus a matrix factorization argument from \cite{bib:methesis},
to ``fill in'' $\calE$ at each of the points of
$\overline{X_1} \setminus
 X_1$.
The quasi-unipotence theorem, conjectured
by Crew \cite[Section~10]{bib:crew2}, follows from the work of any of
Andr\'e \cite{bib:and}, Mebkhout \cite{bib:meb}, or the author
\cite{bib:me7}; see Proposition~\ref{prop:monodromy} for the formulation
we need here.

Theorem~\ref{thm:main} ``almost'' yields a more precise statement 
proposed by
Crew \cite{bib:crew2}, by implying that there exists a finite morphism
$f: X_1 \to X$ such that $f^* \calE$ is unipotent at each point of
$\overline{X_1} \setminus X_1$. (In fact, we will prove Theorem~\ref{thm:main} by
proving this first.) The caveat is that Crew actually wanted
$f$ to be \'etale. We have only been able to achieve this
in the unit-root case, and do not know whether it should be
possible in general; it definitely fails if we allow
$\calE$ to be a log-$F^a$-isocrystal to start with. We discuss these
matters in Section~\ref{sec:etale}.

We hope to extend Theorem~\ref{thm:main} to higher
dimensional varieties in subsequent work. 
In this case, we must allow the morphism $f: X_1 \to X$
to be generically \'etale, but not necessarily finite, because the
normalization of $X$ in a finite extension of its function field
need not be smooth.
For this and other reasons, the argument in higher dimensions
will be technically more involved. 

The precise higher-dimensional analogue of Theorem~\ref{thm:main}
has been formulated conjecturally by Shiho
\cite[Conjecture~3.1.8]{bib:shiho};
by Shiho's work, proving this statement would
provide good comparison results between rigid and crystalline cohomology.
For example, it would yield an alternate proof of finite dimensionality
of rigid cohomology of a curve with coefficients in an overconvergent
$F$-isocrystal, via comparison with crystalline cohomology;
Theorem~\ref{thm:main} does this for curves.
(Crew's proof in \cite{bib:crew2} uses $p$-adic
functional analytic techniques; the general finiteness proof in
\cite{bib:me8} uses a devissage to the curve case.)

\subsection*{Acknowledgments}

Thanks to Pierre Berthelot,
Johan de Jong and Fabien Trihan
for helpful discussions. The author was partially
supported by a National Science Foundation Postdoctoral Fellowship.

\section{Definitions and notations}


We set up notation following \cite{bib:me8}; this makes it a bit more
convenient to work with global objects than does the notation of
\cite{bib:me7}. We retain the convention of \cite{bib:me7} that all matrices
are $n \times n$ matrices and $I$ denotes the identity matrix.

Let $k$ be a field of characteristic $p>0$, and let $C(k)$ be a Cohen
ring for $k$, that is, a complete discrete valuation ring with 
residue field $k$ and maximal ideal generated by $p$.
(See \cite{bib:bou} for proof of existence and basic properties of
Cohen rings.)
Let $K$ be a finite totally ramified extension of $\Frac C(k)$, let
$\calO$ be the integral closure of $C(k)$ in $K$,
and let $v_p$ denote the $p$-adic valuation on $K$.
Assume that there exists a ring
endomorphism $\sigma_0$ on $\calO$ lifting the $p$-power Frobenius on $k$.
Let $q = p^a$ be a power of $p$, and put $\sigma = \sigma_0^a$.

The ring $\calR_r$ 
consists of bidirectional power series $\sum_{i \in \ZZ} c_i u^i$,
with $c_i \in K$, such that
\[
\lim_{i \to \pm \infty} s v_p(c_i) + i  = \infty
 \qquad ( 0 < s \leq r);
\]
for each $s$, the function $w_s(\sum_i c_i u^i) = \min_i \{s v_p(c_i) + i\}$
is a nonarchimedean valuation on $\calR_r$.
The ring $\calR$ (the ``Robba ring'') is the union of the $\calR_r$
over all $r>0$.
Its subring $\calR^{\inte}$ 
consists of those series with $c_i \in \calO$ for all $i$;
this subring is a (noncomplete) discrete valuation ring, unramified
over $\calO$, with residue field the field $k((t))$ of formal Laurent
series in $k$. (In \cite{bib:me7}, the rings $\calR$ and $\calR^{\inte}$
are called $\Gancon$ and $\Gcon$, respectively.) By adding the superscript
$+$ or $-$ to $\calR$ or $\calR^{\inte}$,
we will mean the subring with only nonnegative or nonpositive powers of $u$,
respectively.

For $L/k((t))$ finite separable, there is a natural
discrete valuation ring $\calR^{\inte}_L$,
integral and unramified over $\calR^{\inte}$, with residue field $L$. Namely,
take any
monic polynomial $P(x)$ over $\calR^{\inte}$ 
whose reduction $\overline{P}$
satisfies $K \cong k((t))[x]/(\overline{P}(x))$, and put
$\calR^{\inte}_L \cong \calR^{\inte}[x]/(P(x))$.

The Monsky-Washnitzer algebra of rank $n$ is defined as
\[
W_n = \left\{ \sum_I c_I x^I: c_I \in \calO, \quad \liminf_{I}
\frac{v_p(c_I)}{\sum I} > 0 \right\},
\]
where $I = (i_1, \dots, i_n)$ represents an $n$-tuple of nonnegative integers,
$x^I = x_1^{i_1} \cdots x_n^{i_n}$ and $\sum I = i_1 + \cdots + i_n$.
An \emph{integral dagger algebra} 
is any quotient $A^{\inte}$ of
a Monsky-Washnitzer algebra which is flat over $\calO$ and for which
$\Spec (A^{\inte} \otimes_{\calO} k)$ is smooth over $\Spec(k)$.
(Given $A^{\inte} \otimes_{\calO} k$, one can always find a corresponding
$A^{\inte}$; see \cite{bib:vdp}.)
A \emph{dagger algebra} $A$ is an algebra of the form $A^{\inte}
\otimes_\calO K$ for some integral dagger algebra $A^{\inte}$
(uniquely determined by $A$ and the $p$-adic valuation on $A$).

Given a dagger algebra $A$ with $A^{\inte} \cong W_n/\gotha$
 and $f \in A$ not a zero divisor, the
\emph{localization} $A'$ of $A$ at $f$ is the dagger algebra
with
\[
(A')^{\inte} \cong W_{n+1}/(W_{n+1}\gotha + (fx_{n+1} - 1)W_{n+1});
\]
this is a dagger algebra in which $f$ is invertible.

Given a dagger algebra $A$ with $A^{\inte} \cong W_n/\gotha$ for some
ideal $\gotha$ of $W_n$, we define $\Omega^1_{A/K}$ as the free
$A$-module generated by symbols $dx_1, \dots, dx_n$, modulo relations
of the form $da = 0$ for $a \in \gotha \otimes_{\calO} K$. By construction,
$\Omega^1_{A/K}$ is equipped with a $K$-linear derivation 
$d: A \to \Omega^1_{A/K}$.

\section{Log-$F$-isocrystals and log-$(\sigma, \nabla)$-modules}

Let $A$ be a dagger algebra and 
$\sigma_0: A^{\inte} \to A^{\inte}$ a lift of the
$p$-power Frobenius map extending the given $\sigma_0$ on $\calO$;
again, set $\sigma = \sigma_0^a$.
(Such a lift always exists: again, see \cite{bib:vdp}.)
Given $u \in A$ such that
$u^\sigma/u^q$ is invertible in $A$, we define the logarithmic
module of differentials $\Omega^1_{A/K}[d\log u]$ by adding to
$\Omega^1_{A/K}$ a symbol $d\log u$ such that $u(d\log u) = du$; then
$d\sigma$ extends to $\Omega^1_{A/K}[d \log u]$ sending
$du/u$ to $q\,du/u + d(u^\sigma/u^q)/(u^\sigma/u^q)$.
We define a
\emph{log-$(\sigma, \nabla)$-module} over $A$ (with respect to $u$)
as a finite locally free $A$-module
$M$ equipped with a $\sigma$-linear map $F$ that induces an isomorphism
$F: M \otimes_{A, \sigma} A \to M$,
and with an $A$-linear connection $\nabla: M \to M \otimes_{A}
\Omega^1_{A/K}[d\log u]$ which is integrable (i.e., which satisfies
$\nabla_1 \circ \nabla = 0$ for $\nabla_1: M \otimes \Omega^1_{A/K}[d\log u]
\to M \otimes \wedge^2_A \Omega^1_{A/K}[d\log u]$ induced by $\nabla$)
and which makes the following 
diagram commute:
\[
\xymatrix{
M \ar^-{\nabla}[r] \ar^{F}[d] & M \otimes \Omega^1_{A/K}[d\log u]
\ar^{F \otimes d\sigma}[d] \\
M  \ar^-{\nabla}[r] & M \otimes \Omega^1_{A/K}[d\log u]
}
\]
If $u = 1$, we drop the ``log'' and simply refer to $M$ as
a $(\sigma, \nabla)$-module.
We analogously define $(\sigma, \nabla)$-modules over $\calR$ for $\sigma:
\calR \to \calR$ induced by the $a$-th composition power of a map
$\sigma_0: \calR^{\inte} \to \calR^{\inte}$ of the form
\[
\sum_i c_i u^i \mapsto \sum_i c_i^{\sigma_0} (u^{\sigma_0})^i
\]
lifting the $p$-th power map.
In this case, we take $\Omega^1_{\calR/K}$ to be the free $\calR$-module
generated by $du$. (Note: since we are only considering curves,
the integrability condition will be superfluous in our situations,
as $\wedge^2 \Omega^1$ will vanish.)

Log-$(\sigma, \nabla)$-modules are Zariski-local
avatars of more global objects, namely overconvergent
log-$F^a$-isocrystals. Rather than take space for a full-blown
definition of log-$F^a$-isocrystals here, we summarize the key features
of the definition below.
\begin{itemize}
\item If $X$ is smooth over $k$
and equipped with the fine log structure associated
to some strict normal crossings divisor $Z$, there is a category of
overconvergent
log-$F^a$-isocrystals on the pair $(X,Z)$. (We drop the ``log'' if $Z$
is empty.)
\item An oveconvergent
 log-$F^a$-isocrystal on $(X,Z)$ can be specified by giving
overconvergent
log-$F^a$-isocrystals on an affine cover of $X$ plus isomorphisms on
the pairwise intersections satisfying the
cocycle condition. (Loosely put, the category is a Zariski sheaf.)
\item Given a dagger algebra $A$, a Frobenius lift $\sigma$,
and an element $u \in A^{\inte}$ such that
$u^\sigma/u^q$ is invertible in $A$, the category of 
log-$(\sigma, \nabla)$-modules with respect to $u$
is canonically equivalent to the category
of overconvergent
log-$F^a$-isocrystals on $(X,Z)$, where $X = \Spec (A^{\inte} \otimes_{\calO}
k)$
and $Z \subseteq X$ is the zero locus of $u$. In particular, the former
does not depend on $\sigma$; in fact, there is an explicit formula for
transforming a Frobenius structure with respect to a given $\sigma$ into
a Frobenius structure with respect to another (with respect to the same
$\nabla$).
\end{itemize}
See \cite[Section~6]{bib:ci} for an informal overview of log-$F^a$-isocrystals;
for a much more detailed study, see \cite[Chapter~2]{bib:shiho}.
(It might help to consider the situation without logarithmic
structures first; see \cite{bib:ber} for an introduction there.)

The main input into this paper is the $p$-adic local monodromy theorem
(``Crew's conjecture''), established separately by Andr\'e \cite{bib:and},
Mebkhout \cite{bib:meb}, and the author \cite{bib:me7}.
We say an extension of $k((t))$ is \emph{nearly finite separable} if it is
finite separable over $k^{1/p^n}((t))$ for some nonnegative integer $n$.
With that definition, the local monodromy theorem (e.g., in the form of
\cite[Theorem~6.12]{bib:me7}) implies the following.
\begin{prop} \label{prop:monodromy}
Let $M$ be a $(\sigma, \nabla)$-module over $\calR$.
Then there exists a nearly finite separable extension $L/k((t))$
so that $M \otimes_{\calR^{\inte}} \calR^{\inte}_L$
admits a basis $\bv_1, \dots, \bv_n$ such that
$\nabla \bv_i \in \Span(\bv_1, \dots, \bv_{i-1}) \otimes \Omega^1$
for $i=1,\dots,n$.
\end{prop}
In fact, if $u \in \calR^{\inte}_L$ lifts a uniformizer of $L$ and $K'$
is the integral closure of $K$ in $\calR^{\inte}_L$ (so that
$\calR^{\inte}_L$ is isomorphic to the integral Robba ring with
coefficient field $K'$), then
one can ensure that in fact $\nabla \bv_i \in (K' \bv_1 + \cdots
+ K' \bv_{i-1}) \otimes \frac{du}{u}$; this implies that
$F\bv_i \in K' \bv_1 + \cdots + K' \bv_n$ for all $i$.
If such a basis already exists
over $\calR$, we say $M$ is \emph{unipotent} over $\calR$.

\section{Matrix factorizations}

In this section, let $A$ be a dagger algebra such that
$\Spec (A \otimes_{\calO} k)$ is a smooth affine geometrically connected
curve over $\Spec(k)$; in particular, $A$ is an integral domain. 
Suppose $t \in A \otimes_{\calO}
k$ generates a prime ideal with residue field $k$;
choose a lift $u \in A^{\inte}$
of $t$. Then there is a natural embedding $\rho_u: A^{\inte} \hookrightarrow
\calR^{+,\inte}$
sending $u$ to the series parameter; this extends to an embedding of any
localization of $A$ into $\calR$. We identify each localization of $A$
with its image under $\rho_u$.

\begin{lemma} \label{lem:approx}
For $r>0$, let $U$ be an invertible matrix over $\calR_r$.
Then there exists an
invertible matrix $V$ over some localization $A'$ of $A$ such that
$w_r(VU - I) > 0$.
\end{lemma}
\begin{proof}
Let $A_1$ be the localization of $A$ at $u$; then $A_1$
is dense in $\calR_r$ under $w_r$.
Thus we can choose a matrix
 $V$ over $A_1$ such that
$w_r(V - U^{-1}) > -w_r(U)$; then 
\[
w_r(VU-I) = w_r((V-U^{-1})U) \geq w_r(V - U^{-1}) + w_r(U) > 0.
\]
Since $w_r(VU - I) > 0$, we have $w_r(\det(VU) - 1) > 0$, so in
particular $\det(VU) \neq 0$. Hence $\det(V) \neq 0$, so we can form
the localization $A'$ of $A_1$ at $\det(V)$. Over
$A'$, $V$ becomes an invertible matrix, as desired.
\end{proof}

\begin{prop} \label{prop:matfact}
Let $U$ be an invertible matrix over $\calR$.
Then there exist
invertible matrices $V$ over some localization
$A'$ of $A$ and $W$ over $\calR^+$
such that $U = VW$.
\end{prop}
\begin{proof}
Choose $r>0$ so that $U$ is invertible over $\calR_r$.
By Lemma~\ref{lem:approx}, we can find an invertible
 matrix $X$ over some localization
$A_1$ of $A$ such that $w_r(XU - I) > 0$. By
\cite[Proposition~6.5]{bib:me7}, we can write $XU$ as a product $YZ$
with $Y$ invertible over $\calR^-$
and $Z$ invertible over $\calR^+$.
Let $A'$ be the localization of $A_1$ at $u$;
then $\calR^- \subseteq
A'$, so $X^{-1}Y$ is invertible over $A'$. Thus we may take
$V = X^{-1}Y$ and $W = Z$.
\end{proof}

\section{Semistable reduction}

In this section, we prove that an overconvergent
$F^a$-isocrystal $\calE$ on $X$ which is
unipotent at each point of $\overline{X}-X$, for $\overline{X}$ a smooth
compactification of $X$, has ``semistable reduction''. This will yield
our proof of Theorem~\ref{thm:main}.

We begin with a result that translates unipotence of a $(\sigma,
\nabla)$-module
into semistable reduction. The argument is based on the proof of
\cite[Theorem~5.0.1]{bib:methesis}.
\begin{theorem} \label{thm:semistab}
Let $A$ be a dagger algebra equipped with a Frobenius lift $\sigma$.
Suppose the image of $u \in A^{\inte}$
in $A^{\inte} \otimes_{\calO} k$ generates a prime ideal and
$u^\sigma/u^q$ is a unit in $A$. Let $A'$ be the localization
of $A$ at $u$, and let $M$
be a free $(\sigma, \nabla)$-module over $A'$ which becomes unipotent
over $\calR$ (where $A'$ is identified with a subring of $\calR$
via $\rho_u$).
Then $M$ is isomorphic to a log-$(\sigma, \nabla)$-module,
with respect to $u$, over some localization $A''$ of $A$ in which 
$u$ is not invertible.
\end{theorem}
\begin{proof}
Let $\be_1, \dots, \be_n$ be a basis of $M$, so that
$F\be_j = \sum_i \Phi_{ij} \be_i$ and $\nabla \be_j = \sum_i N_{ij} \be_i
\otimes \frac{du}{u}$.
Define the differential operator $\theta(f) = u \frac{df}{du}$.
By hypothesis, there exists a matrix $U$ over $\calR$ such that
$U^{-1} \Phi U^\sigma$ and $U^{-1} N U + U^{-1} \theta(U)$ 
have entries in
$\calO$.
By Proposition~\ref{prop:matfact}, we can factor $U$ as $VW$,
where $V$ is invertible over
some localization $A_1$ of $A'$
and $W$ is invertible over $\calR^+$.
Now put $\bv_j = \sum_{j} V_{ij} \be_i$; then
$F\bv_j = \sum_i \tilde{\Phi}_{ij} \bv_i$
and $\nabla \bv_j = \sum_i \tilde{N}_{ij} \bv_i
\otimes \frac{du}{u}$, where
\begin{align*}
\tilde{\Phi} &= V^{-1} \Phi V^\sigma = W (U^{-1} \Phi U^\sigma) W^{-\sigma} \\
\tilde{N} &= V^{-1} N V + V^{-1} \theta(V) = 
W(U^{-1}NU+U^{-1} \theta(U))W^{-1} - \theta(W) W^{-1}
\end{align*}
have entries in $A'' = A_1 \cap \calR^+$, which is a localization of $A$
(because it contains $A$ and is contained in the localization $A_1$)
in which $u$ is not invertible (because $u$ is not invertible in
$\calR^+$).
\end{proof}

We now proceed to the proof of our main theorem.
\begin{proof}[Proof of Theorem~\ref{thm:main}]
Let $K(X)$ be the function field of $X$, and let $\overline{X}$ be a smooth
compactification of $X$. 
Without loss of generality, enlarge $k$ so that
the geometric points of $Z = \overline{X} 
\setminus X$ are $k$-rational. 
(Given the desired result over a finite extension of $k$,
we deduce the result over $k$ by restriction of scalars.)
For each geometric point $x$ of $Z$, choose a function $t_x \in K(X)$ with
a simple zero at $x$, and choose an open affine neighborhood
$U_x$ of $x$ in $\overline{X}$ such that $U_x \cap Z = \{x\}$
and $\divi(t_x) \cap U_x = \{x\}$. Let $A_x$ be a dagger algebra with
$U \cong \Spec(A^{\inte}_x 
\otimes_{\calO} k)$, choose a Frobenius lift $\sigma$ on
$A_x$ and choose a lift $u_x$ of $t_x$
in $A_x^{\inte}$. Let $A'_x$ be the localization of $A_x$ at $u_x$;
then the restriction of
$\calE$ to $U \setminus \{x\}$ corresponds to
a $(\sigma, \nabla)$-module $M_x$ over $A'_x$.
After shrinking $U_x$ if needed, we may assume that $M_x$
is free over $A'_x$.

Let $\rho_x$ be the embedding of $A_x$ into $\calR^{+,\inte}$ 
sending $u_x$ to the series parameter. By Proposition~\ref{prop:monodromy},
there exists a nearly finite separable
extension $L_x$ of the $t_x$-adic completion of $K(X)$
such that $M_x \otimes_{A'_x} \calR^{\inte}_{L_x}$ is unipotent.
By Krasner's lemma, after replacing $k$ with $k^{1/p^n}$ for some nonnegative
integer $n$, we can choose a finite separable
extension $L$ of $K(X)$ whose completion at any
point above $x$ contains $L_x$ for each $x \in Z$.
After enlarging $k$ again, we may assume that the places of $L$ above $x$
are all $k$-rational.
Let $f: \overline{X_1} \to \overline{X}$ 
be the cover corresponding to the extension $L/K(X)$,
and put $X_1 = f^{-1}(X)$.

For each geometric point $y$ of $Z_1 = \overline{X_1} \setminus X_1$,
choose a function $t_y \in K(X_1)$ with
a simple zero at $y$ and no multiple zeroes; then $t_y$ gives rise to
a map $g_y: \overline{X_1}
 \to \PP^1$. (If $k$ is finite, it may be necessary to
enlarge it again to find such $t_y$.)
Put $x = f(y)$,
and choose an open affine neighborhood
$V_y$ of $y$ in $\overline{X_1}$ such that $V_y \cap Z_1 = \{y\}$,
$\divi(t_y) \cap V_y = \{y\}$, $f(V_y) \subseteq U_x$,
and $V_y$ does not meet the branch locus
of $g_y$.
Then there is a dagger algebra $B_y$ with
$V_y \cong \Spec(B^{\inte}_y \otimes_{\calO} k)$ which is
a localization of a finite extension of $A_x$.

Choose a lift $u_y$ of $t_y$ in $B_y^{\inte}$. Since $V_y$ is unramified
over its image under $g_y$, $B_y$ is
finite and unramified over some localization $C_y$ of its
subring $K \langle u_y \rangle^\dagger$ (the $p$-adic closure of
$K[u_y]$ within $B_y$). Now $K \langle u_y \rangle^\dagger$ admits
a $p$-power Frobenius lift $\sigma'_0$ sending $u_y$ to $u_y^\sigma$;
this lift extends to the localization $C_y$, then to the unramified
extension $B_y$. Put $\sigma' = (\sigma'_0)^a$.

Let $B'_y$ be the localization of $B_y$ at $u_x \in A_x \subseteq B_y$,
which is the same as the localization at $u_y$ because 
$\divi(t_x) \cap V_y = \divi(t_y) \cap V_y = \{y\}$,
and put $V'_y = \Spec((B'_y)^{\inte} \otimes_{\calO} k) =
 V_y \setminus \{y\}$;
then the restriction of $f^* \calE$ to $V'_y$
corresponds to the $(\sigma', \nabla)$-module
$M_x \otimes_{A'_x} B'_y$. (That is, its connection is the one induced
from $M_x$, but its Frobenius structure is defined with respect to
$\sigma'$ instead of $\sigma$.) By construction, this $(\sigma',
\nabla)$-module is unipotent. Moreover, $u_y^{\sigma'}/u_y^q = 1$,
so Theorem~\ref{thm:semistab} implies that $M_x \otimes_{A'_x}
B'_y$ is isomorphic to a log-$(\sigma, \nabla)$-module over
some localization $B''_y$ of $B_y$ in which $u_y$ is not invertible.
If $V''_y = \Spec((B''_y)^{\inte} \otimes_{\calO} k)$, then
$V''_y$ is an open affine neighborhood of $y$ in $\overline{X_1}$
on which $f^*\calE$ extends to a log-$F$-isocrystal.

In short, we have an open affine neighborhood of each $y \in Z_1$
in $\overline{X_1}$,
on which $f^* \calE$ extends to a log-$F$-isocrystal relative
to $\{y\}$. Each
neighborhood contains no other points of $y$, so the pairwise
intersections all lie in $X_1$. Thus we automatically have glueing
isomorphisms on the log-$F$-isocrystals satisfying the cocycle
conditions (since $f^* \calE$ is defined on $X_1$), yielding
a log-$F$-isocrystal on $(\overline{X_1}, Z_1)$, as desired.
\end{proof}

\section{Finite versus \'etale}
\label{sec:etale}

As noted earlier, Crew \cite[Section~10]{bib:crew2}
conjectured that an overconvergent $F$-isocrystal on a curve
should extend to a log-$F$-isocrystal after a base extension
which is not just finite and generically \'etale, but actually \'etale.
It is unclear whether this should hold in general; the
best we can do at the moment is prove it in the unit-root case,
as done below. Note that
this proof does not use the full strength of the quasi-unipotence theorem,
but only the unit-root case; this case is due to Tsuzuki \cite{bib:tsu1}.
Also note that a unit-root log-$F$-isocrystal is automatically an
$F$-isocrystal, so there is no ``log'' in the statement of the theorem.

The proof of the following lemma is straightforward.
\begin{lemma} \label{lem:unr}
Let $B$ be a matrix over $k [[ t ]]$, for $k$
a perfect field of characteristic $p>0$, and
let $\tau$ denote the $q$-th power map. Then
any solution $D$ of either of the matrix equations
\[
D^{-1} B D^\tau = I \qquad \mbox{or} \qquad D^\tau - D = B
\]
over the integral closure of $k[[t]]$ in 
$k((t))^{\alg}$ is defined over an unramified extension
of $k[[t]]$.
\end{lemma}

\begin{theorem}
Let $X$ be a smooth, geometrically connected curve over a perfect field $k$ 
of characteristic $p>0$, and let $\calE$ be an
overconvergent unit-root $F^a$-isocrystal on $X/K$. Then there exists
a finite \'etale morphism $f: X_1 \to X$, a smooth compactification
$j: X_1 \hookrightarrow \overline{X_1}$ of $X_1$, and
a unit-root $F^a$-isocrystal $\calF$ on $\overline{X_1}$
such that $j^*\calF \cong f^* \calE$.
\end{theorem}
\begin{proof}
If $X$ is projective, there is nothing to prove, so we assume $X$ is affine.
Let $A$ be a dagger algebra with $X \cong \Spec(A^{\inte} \otimes_\calO k)$,
and choose a Frobenius lift $\sigma$ on $A$. Then $\calE$ corresponds
to a $(\sigma, \nabla)$-module $M$ over $A$. Choose (not necessarily free) generators $\bv_1, \dots, \bv_m$ of $M$,
and let $N \subset M$ be the $A^{\inte}$-span of $F^i \bv_j$ over
$i=0,1,\dots$ and $j=1,\dots, m$. Then $N$ is locally free over
$A^{\inte}$.

Let $L$ be the $p$-adic
completion of the valuation subring of $\Frac A$;
note that $N$ is free over $L$
and $F$ acts on any basis of $N$ over $L$ via an invertible matrix.
Let $\pi$ be a uniformizer of $\calO$,
and pick an integer $d$ such that $v_p(\pi^d) > 1/(p-1)$.

Given any basis
$\be_1, \dots, \be_n$ of $N$ to start with,
define the matrix $\Phi$ by $F\be_j = \sum_{ij} \Phi_{ij} \be_i$. 
We then solve the matrix equation
$C_i^{-1} \Phi C_i^\sigma \equiv I \pmod{\pi^i}$ for $i=1, \dots, d$
to obtain a matrix $C_d$ over some finite unramified
extension $L'$ of $L$. Then $N \otimes_{A^{\inte}} L'$ 
admits a basis $\bw_1, \dots, \bw_n$ for which
$F\bw_i \equiv \bw_i \pmod{\pi^d}$.

If we insist that $L'$ be minimal for the existence of the basis of the
desired form, then it is unique; 
in particular, it does not depend on the choice of the
starting basis. Let $X_1$ be a curve for which $K(X_1) \cong L'/\pi L'$
and let $f: X_1 \to X$ be the induced map.
If we choose the initial basis $\be_1, \dots, \be_n$
over a localization $A^{\inte}_1$ of $A^{\inte}$
over which $N$ becomes free, we discover that $f$ is \'etale over
any point in $\Spec(A_1^{\inte} \otimes_{\calO} k)$ by
Lemma~\ref{lem:unr}. (Namely, $C_0$ satisfies an equation modulo $\pi$
of the
first type in the lemma, while
 $C_{i-1}^{-1} C_i = I + \pi^{i-1} D$ for some matrix
$D$ satisfying an equation modulo $\pi$ of the second type.) Since $N$ is
locally free over $A^{\inte}$,
we can arrange for $\Spec(A_1^{\inte} \otimes_{\calO} k)$
to contain any closed point of $X$. Thus $f$ is finite \'etale.
Since $v_p(\pi^d) > 1/(p-1)$, we may apply \cite[Theorem~5.1.1]{bib:tsu1}
(at least for $a=1$; see \cite[Proposition~6.11]{bib:me7} for
a reduction to this case) to see that $f^* \calE$ admits
a basis of elements in the kernel of $\nabla$, on which $F$ acts
by a matrix over $\calO$. This allows us to extend $f^* \calE$
to a compactification $\overline{X_1}$ of $X_1$, as desired.
\end{proof}

The general case is complicated by the fact that there is no obvious
global criterion for when a non-unit-root overconvergent $F$-isocrystal
becomes unipotent.
Moreover, a slightly
stronger statement turns out to be false: it is not always possible to pull
back a log-$F$-isocrystal by an \'etale map to obtain a log-$F$-isocrystal
that extends over a smooth compactification. This is demonstrated by an 
example of Tsuzuki \cite[Section~6]{bib:tsu3}. In this example
(the ``Bessel isocrystal''), an overconvergent
$F$-isocrystal over $\PP^1-\{0, \infty\}$ is unipotent over $0$ and
becomes unipotent over $\infty$ after making a tamely ramified quadratic
extension. Thus the $F$-isocrystal already extends to a log-$F$-isocrystal
over $\PP^1-\{\infty\}$ with respect to the divisor $(0)$, but no
\'etale cover of $\PP^1-\{\infty\}$ can induce anything but a totally
wildly ramified extension of the local ring at $\infty$.


\begin{thebibliography}{Ke2}
\bibitem[A]{bib:and}
Y. Andr\'e, Filtrations de type Hasse-Arf et monodromie $p$-adique, 
\textit{Invent. Math.} \textbf{148} (2002), 285--317.

\bibitem[Be]{bib:ber}
P. Berthelot, G\'eom\'etrie rigide et cohomologie des vari\'et\'es
algebriques de caract\'eristique $p$, Introductions aux cohomologies
$p$-adiques (Luminy, 1984), \textit{M\'em. Soc. Math. France} \textbf{23}
(1986), 7--32.

\bibitem[Bo]{bib:bou}
N. Bourbaki, \textit{Algebre Commutative, chap.\ IX--X}, Masson (Paris), 1983.

\bibitem[CI]{bib:ci}
R. Coleman and A. Iovita, Hidden structures on semistable curves,
preprint available at \texttt{math.berkeley.edu/\~{}coleman}.

\bibitem[Cr]{bib:crew2}
R. Crew, Finiteness theorems for the cohomology of an overconvergent isocrystal
on a curve, \textit{Ann. Scient. \'Ec. Norm. Sup.} \textbf{31} (1998), 717--763.


\bibitem[dJ]{bib:dej2}
A.J. de~Jong, Barsotti-Tate groups and crystals,
Proceedings of the International Congress of Mathematicians,
\textit{Doc. Math.} Extra Vol.\ II (1998),
259--265.

\bibitem[Ke1]{bib:methesis}
K.S. Kedlaya, Descent theorems for overconvergent $F$-crystals,
Ph.D. thesis, Massachusetts Institute of Technology, 2000.

\bibitem[Ke2]{bib:me7}
K.S. Kedlaya, A $p$-adic local monodromy theorem, to appear in
\textit{Ann. Math.}; \texttt{arXiv: math.AG/0110124}.

\bibitem[Ke3]{bib:me8}
K.S. Kedlaya, Finiteness of rigid cohomology with coefficients, preprint;
\texttt{arXiv: math.AG/0208027}.


\bibitem[M]{bib:meb}
Z. Mebkhout, Analogue $p$-adique du Th\'eor\`eme de Turrittin
et le Th\'eor\`eme de la monodromie $p$-adique, 
\textit{Invent. Math.} \textbf{148} (2002), 319--351.

\bibitem[S]{bib:shiho}
A. Shiho, Crystalline fundamental groups. II. Log convergent cohomology
and rigid cohomology, \textit{J. Math. Sci. Univ. Tokyo} \textbf{9} (2002),
1--163.


\bibitem[T1]{bib:tsu1}
N. Tsuzuki, Finite local monodromy of
overconvergent unit-root $F$-crystals on a curve,
\textit{Amer. J. Math.} \textbf{120} (1998), 1165--1190.

\bibitem[T2]{bib:tsu3}
N. Tsuzuki, Slope filtration of quasi-unipotent overconvergent
$F$-isocrystals, \textit{Ann. Inst. Fourier, Grenoble} \textbf{48} (1998),
379--412.


\bibitem[vdP]{bib:vdp}
M. van der Put, The cohomology of Monsky and Washnitzer, Introductions
aux cohomologies $p$-adiques (Luminy, 1984),
\textit{M\'em. Soc. Math. France} \textbf{23} (1986), 33--60.

\end{thebibliography}
\end{document}